\documentclass{article}

\usepackage{color}

\usepackage{amsmath}
\usepackage{amssymb}
\usepackage{amsthm}
\usepackage{amsfonts}
\usepackage{latexsym}
\usepackage{amssymb}
\usepackage{url}
\usepackage[dvipdfmx]{graphicx}

\DeclareMathOperator*{\grad}{grad}
\newcommand{\real}{\mathbb{R}}

\newcommand{\norm}[1]{\left\lVert{#1}\right\rVert}
\newcommand{\ip}[1]{\left\langle{#1}\right\rangle}

\newtheorem{definition}{Definition}[section]
\newtheorem{theorem}{Theorem}[section]

\newtheorem{problem}{Problem}[section]

\begin{document}

\title{Hybrid Riemannian Conjugate Gradient Methods with Global Convergence Properties
\thanks{This work was supported by JSPS KAKENHI Grant Number JP18K11184.}
}
\author{Hiroyuki Sakai \and Hideaki Iiduka}
\maketitle

\begin{abstract}
This paper presents Riemannian conjugate gradient methods and global convergence analyses under the strong Wolfe conditions.
The main idea of the proposed methods is to combine the good global convergence properties of the Dai-Yuan method with the efficient numerical performance of the Hestenes-Stiefel method.
One of the proposed algorithms is a generalization to Riemannian manifolds of the hybrid conjugate gradient method of the Dai and Yuan in Euclidean space.
The proposed methods are compared well numerically with the existing methods for solving several Riemannian optimization problems.
\end{abstract}

\section{Introduction}
This paper focuses on the conjugate gradient method. Nonlinear conjugate gradient methods in Euclidean space are a class of important methods for solving unconstrained optimization problems. In \cite{hestenes1952}, Hestenes and Stiefel developed a conjugate gradient method for solving linear systems with a symmetric positive-definite matrix of coefficients. In \cite{fletcher1964}, Fletcher and Reeves extended the conjugate gradient method to unconstrained nonlinear optimization problems. Theirs is the first nonlinear conjugate gradient method in Euclidean space. Al-Baali \cite{baali1985} indicated that the Fletcher-Reeves method converges globally and generates the descent direction with an inexact line search when the step size satisfies the strong Wolfe conditions \cite{wolfe1969,wolfe1971}. Polak and Ribi\`ere \cite{polak1969} introduced a conjugate gradient method with good numerical performance. Dai and Yuan \cite{dai1999} introduced a conjugate gradient method with a better global convergence property than that of the Fletcher-Reeves method. The Hestenes-Stiefel and Polak-Ribi\`ere-Polyak methods do not always converge under the strong Wolfe conditions, and for this reason, hybrid conjugate gradient methods have been presented in \cite{dai2001,hu1991,touati1990}. Touati-Ahmed and Storey \cite{touati1990}, and Hu and Storey \cite{hu1991} proposed methods combining the Fletcher-Reeves and Polak-Ribi\`ere-Polyak methods. Moreover, Dai and Yuan \cite{dai2001} proposed the hybrid conjugate gradient method, which combines the Dai-Yuan method and the Hestenes-Stiefel method. These nonlinear conjugate gradient methods in Euclidean space are summarized by Hager and Zhang in \cite{hager2006}.

The conjugate gradient method in Euclidean space is applicable to a Riemannian manifold. In \cite{smith1994}, Smith introduced the notion of Riemannian optimization using the exponential map and parallel translation. However, using the exponential map or parallel translation on a Riemannian manifold is generally not computationally efficient. Absil, Mahony, and Sepulchre \cite{absil2008} proposed to use a mapping called a retraction that approximates the exponential map. Moreover, they introduced the notion of vector transport, which approximates parallel transport. In addition, Ring and Wirth \cite{ring2012} introduced generalized line search methods (e.g., the Wolfe conditions \cite{wolfe1969,wolfe1971}) on Riemannian manifolds.

Using the retraction and vector transport, Ring and Wirth \cite{ring2012} presented a Fletcher-Reeves type of nonlinear conjugate gradient method on Riemannian manifolds. They indicated that the Fletcher-Reeves methods have a global convergence property under the strong Wolfe conditions. However, their convergence analysis assumed that the vector transport does not increase the norm of the search direction vector, which is not the standard assumption (see \cite[Section 5]{sato2013}). To remove this unnatural assumption, Sato and Iwai \cite{sato2013} introduced the notion of scaled vector transport \cite[Definition 2.2]{sato2013}. They proved that by using scaled vector transport, the Fletcher-Reeves method on a Riemannian manifold generates a descent direction at every iteration and converges globally without impractical assumptions. Similarly, in \cite{sato2015}, Sato used scaled vector transport in a convergence analysis. He indicated that the Dai-Yuan-type Riemannian conjugate gradient method generates a descent direction at every iteration and converges globally under the Wolfe conditions. This means that the Dai-Yuan method has a better global convergence property than that of the Fletcher-Reeves method on Riemannian manifolds, since the latter has to resort to the \emph{strong} Wolfe conditions, whereas the former only requires the Wolfe conditions.

In this paper, we propose hybrid Riemannian conjugate gradient methods exploiting the idea used in the paper \cite{dai2001}.
One of the methods we propose has already been used in numerical experiments (e.g., \cite[(43)]{hawe2013analysis}, \cite[Table 1]{selvan2012descent}), but no convergence analysis has yet been presented for it.
Our methods combine the good numerical performance of the Hestenes-Stiefel method with the efficient global convergence property of the Dai-Yuan method. Moreover, we present convergence analyses of our methods. The proofs are along the lines of \cite[Theorem 2.3]{dai2001}, except that the step-size assumption is stronger than that of the Euclidean case. This is due to the use of scaled vector transport. Our hybrid methods converge globally if the size of the parameter, which is used to determine the search direction, with respect to that of the Dai-Yuan method is in a certain range (Theorem \ref{thm:main}). We provide two examples which satisfy such a condition. In numerical experiments, we show that our hybrid methods outperform the Dai-Yuan and Polak-Ribi\`ere-Polyak methods.

This paper is organized as follows. Section \ref{sec:RCG} reviews the fundamentals of Riemannian geometry and Riemannian optimization. Section \ref{sec:main} proposes the hybrid Riemannian conjugate gradient methods and presents global convergence analyses for them. Section \ref{sec:exp} compares our methods with the existing Riemannian conjugate gradient methods through numerical experiments. Section \ref{sec:conclude} concludes the paper with mention of future work.

\section{Riemannian Conjugate Gradient Methods}~\label{sec:RCG}
Let us start by reviewing the nonlinear conjugate gradient methods in Euclidean space. The search direction $\eta_k$ of the nonlinear conjugate gradient method is determined by $\eta_0=-{\nabla}f(x_0)$ and
 \begin{align}~\label{eq:EuCGdir}
 \eta_{k+1}=-{\nabla}f(x_{k+1})+\beta_{k+1}\eta_{k},
 \end{align}
where $x_0 \in \real^n$, $\beta_0=0$, and $\beta_k$ is a parameter to be suitably defined. Well-known formulas for $\beta_k$ are the Fletcher-Reeves (FR) \cite{fletcher1964}, Dai-Yuan (DY) \cite{dai1999}, Polak-Ribi\`ere-Polyak (PRP) \cite{polak1969}, and Hestenes-Stiefel (HS) \cite{hestenes1952} formulas, given by
\begin{align}
\label{eq:EuFR}
\beta^\textrm{FR}_k &= \frac{\norm{{\nabla}f(x_k)}^2}{\norm{{\nabla}f(x_{k-1})}^2}, \\
\label{eq:EuDY}
\beta^\textrm{DY}_k &= \frac{\norm{{\nabla}f(x_k)}^2}{\eta^\top _{k-1}y_{k-1}}, \\
\label{eq:EuPRP}
\beta^\textrm{PRP}_k &= \frac{{\nabla}f(x_k)^\top y_{k-1}}{\norm{{\nabla}f(x_{k-1})}^2}, \\
\label{eq:EuHS}
\beta^\textrm{HS}_k &= \frac{{\nabla}f(x_k)^\top y_{k-1}}{\eta_{k-1}^\top y_{k-1}},
\end{align}
respectively, where $y_{k-1}={\nabla}f(x_k)-{\nabla}f(x_{k-1})$.

In the Euclidean space setting, a line search optimization algorithm updates the current iterate $x_k$ to the next iterate $x_{k+1}$ with the updating formula, 
\begin{align}~\label{eq:EuCGupd}
x_{k+1}=x_k+\alpha_k\eta_k,
\end{align}
where $\alpha_k>0$ is a positive step size. One often chooses a step size $\alpha_k>0$ to satisfy the Wolfe conditions \cite{wolfe1969,wolfe1971}, namely,
\begin{align}
\label{eq:EuWolfe1}
f(x_k+\alpha_k\eta_k) \leq f(x_k) + c_1\alpha_k{\nabla}f(x_k)^\top\eta_k, \\
\label{eq:EuWolfe2}
{\nabla}f(x_k+\alpha_k\eta_k)^\top \eta_k \geq c_2{\nabla}f(x_k)^\top\eta_k \textcolor{red}{,}
\end{align}
where $0<c_1<c_2<1$. When the step size satisfies the following condition, which is a substitute of \eqref{eq:EuWolfe2}:
\begin{align}
\label{eq:EusWolfe}
|{\nabla}f(x_k+\alpha_k\eta_k)^\top\eta_k| \leq c_2|{\nabla}f(x_k)^\top\eta_k|,
\end{align}
we call \eqref{eq:EuWolfe1} and \eqref{eq:EusWolfe} the strong Wolfe conditions.

In \cite{dai2001}, Dai and Yuan proved that the method defined by \eqref{eq:EuCGdir} and \eqref{eq:EuCGupd} produces a descent search direction at every iteration and converges globally if the step size $\alpha_k>0$ satisfies \eqref{eq:EuWolfe1} and \eqref{eq:EuWolfe2}, and $\beta_k$ satisfies
\begin{align*}
-\sigma \leq \frac{\beta_k}{\beta_k^\textrm{DY}} \leq 1,
\end{align*}
where $\sigma := (1-c_2)/(1+c_2)$ and $c_2$ is a constant in the second condition \eqref{eq:EuWolfe2}. In this paper, we extend these choices of the parameter $\beta_k$ to Riemannian manifolds.

Now we will briefly outline Riemannian optimization, especially the Riemannian conjugate gradient method, by summarizing \cite{absil2008}. Moreover, we will introduce relevant notation of Riemannian geometry.

Let $(M,g)$ be a Riemannian manifold with a Riemannian metric $g$, and let $T_xM$ be the tangent vector space of $M$ at a point of $x \in M$. In addition, let $TM$ be the tangent bundle of $M$, which is defined by $TM=\bigcup_{x \in M}T_xM$. Let $f:M \rightarrow \mathbb{R}$ be a smooth objective function. Throughout this paper, to simplify the notation, we will write the Riemannian metric $g(\cdot ,\cdot)$ as $\ip{\cdot ,\cdot}$. Given a smooth function $f:M\rightarrow\real$, the gradient of $f$ at a point $x \in M$, denoted by $\grad f(x)$, is defined as the unique element of $T_xM$ that satisfies
\begin{align*}
df_x(\xi)=\ip{\grad f(x),\xi}_x \qquad (\xi \in T_xM).
\end{align*}

An unconstrained optimization problem on a Riemannian manifold $M$ is expressed as follows:
\begin{problem}
Let $f:M\rightarrow\real$ be smooth. Then, we would like to
\begin{align*}
\textrm{minimize}& \quad f(x), \\
\textrm{subject to}& \quad x \in M.
\end{align*}
\end{problem}

In order to generalize line search optimization algorithms to Riemannian manifolds, we will use the notions of a retraction and vector transport (see \cite{absil2008}), which are defined as follows:

\begin{definition}
[Retraction]~\label{def:retraction}
Let $M$ be a manifold and $TM$ be a tangent bundle of a manifold $M$. Any smooth map $R:TM{\rightarrow}M$ is called a retraction on $M$, if it has the following properties.
\begin{itemize}
\item $R_x(0_x)=x$, where $0_x$ denotes the zero element of $T_xM$;
\item With the canonical identification $T_{0_x}T_xM \simeq T_xM$, $R_x$ satisfies $\textrm{D}R_x(0_x)[\xi]=\xi$ for all $\xi \in T_xM$,
\end{itemize}
where $R_x$ denotes the restriction of $R$ to $T_xM$ and $\text{D}R$ is the differential of $R$ (see \cite[Section 3]{absil2008}).
\end{definition}

\begin{definition}
[Vector transport]~\label{def:vector transport}
Let $M$ be a manifold and $TM$ be a tangent bundle of $M$. Any smooth map $\mathcal{T}:TM{\oplus}TM{\rightarrow}TM$, where $\oplus$ denotes the Whitney sum, is called vector transport on $M$, if it has the following properties.
\begin{itemize}
\item There exists a retraction $R$, called the retraction associated with $\mathcal{T}$, such that $\mathcal{T}_{\eta}(\xi) \in T_{R_x(\eta)}M$ for all $x \in M$, and for all $\eta ,\xi \in T_xM$\textcolor{red}{;}
\item $\mathcal{T}_{0_x}(\xi)=\xi$ for all $\xi \in T_xM$;
\item $\mathcal{T}_{\eta}(a\xi + b\zeta)=a\mathcal{T}_{\eta}(\xi)+b\mathcal{T}_{\eta}(\zeta)$ for all $a,b \in \mathbb{R}$, and for all $\eta ,\xi ,\zeta \in T_xM$.
\end{itemize}
where $\mathcal{T}_{\eta}(\xi)$ denotes $\mathcal{T}(\eta ,\xi)$.
\end{definition}

In this paper, we will focus on the differentiated retraction $\mathcal{T}^\textrm{R}$ as a vector transport, defined by
\begin{align}~\label{eq:difft}
\mathcal{T}^{R}_{\eta}(\xi):=\textrm{D}R_{x}(\eta)[\xi] \qquad (\xi \in T_xM)\textcolor{red}{,}
\end{align}
where $x \in M$ and $\eta \in T_xM$. It is easy to prove that $\mathcal{T}^R$ satisfies the properties of Definition \ref{def:retraction} (see \cite[Chapter 8]{absil2008}).

In Riemannian optimization, by using a retraction $R$ and vector transport $\mathcal{T}$ on $M$, we can generalize the updating formula \eqref{eq:EuCGupd} and the search direction of the conjugate gradient method \eqref{eq:EuCGdir} to, respectively, 
\begin{align}
\label{eq:CGupd}
x_{k+1}&=R_{x_k}(\alpha_k\eta_k), \\
\label{eq:CGdir}
\eta_{k+1}&=-{\grad}f(x_{k+1})+\beta_{k+1}\mathcal{T}_{\alpha_{k}\eta_{k}}(\eta_{k}),
\end{align}
where $\alpha_k>0$ is a positive step size (see \cite{absil2008}). We call the search direction $\eta_k$ a descent direction if $\eta_k$ satisfies
\begin{align*}
\ip{{\grad}f(x_k),\eta_k}_{x_k} < 0.
\end{align*}
Moreover, the line search conditions \eqref{eq:EuWolfe1} and \eqref{eq:EuWolfe2} can be generalized to Riemannian manifolds as follows:
\begin{align}
\label{eq:Armijo}
f(R_{x_k}(\alpha_k\eta_k)) \leq f(x_k)+c_1\alpha_k\ip{{{\grad}f(x_k),\eta_k}}_{x_k}, \\
\label{eq:Wolfe}
\ip{ {\grad}f(R_{x_k}(\alpha_k\eta_k)),\textrm{D}R_{x_k}(\alpha_k\eta_k)[\eta_k]}_{R_{x_k}(\alpha_k\eta_k)} \geq c_2\ip{{{\grad}f(x_k),\eta_k}}_{x_k},
\end{align}
where $0<c_1<c_2<1$ (see \cite{sato2015,sato2013}). We call \eqref{eq:Armijo} the Armijo condition. Moreover, the second of the strong Wolfe conditions \eqref{eq:EusWolfe} can be rewritten as
\begin{align}
\label{eq:sWolfe}
\left|\ip{ {\grad}f(R_{x_k}(\alpha_k\eta_k)),\textrm{D}R_{x_k}(\alpha_k\eta_k)[\eta_k]}_{R_{x_k}(\alpha_k\eta_k)}\right| \leq c_2\left|\ip{{{\grad}f(x_k),\eta_k}}_{x_k}\right|.
\end{align}

Sato and Iwai \cite{sato2013} introduced the notion of scaled vector transport. A scaled vector transport of the $k$-th iterate $\mathcal{T}^{(k)}$ associated with $\mathcal{T}^\textrm{R}$ is defined by
\begin{align}~\label{eq:scaled}
\mathcal{T}^{(k)}_{\alpha_k\eta_k}(\eta_k):=
	\begin{cases}
	\mathcal{T}^R_{\alpha_k\eta_k}(\eta_k), & \textrm{if}\, \norm{\mathcal{T}_{\alpha_k\eta_k}^R(\eta_k)}_{x_{k+1}}\leq\norm{\eta_{k}}_{x_k}, \\
	\dfrac{\norm{\eta_k}_{x_k}}{\norm{\mathcal{T}^R_{\alpha_k\eta_k}(\eta_k)}_{R_{\alpha_k\eta_k}(\eta_k)}}\mathcal{T}^R_{\alpha_k\eta_k}(\eta_k), & \textrm{otherwise}.
	\end{cases}
\end{align}
Note that scaled vector transport does not satisfy the properties of Definition \ref{def:vector transport}. Thus, we cannot call this vector transport with mathematical exactitude; however, by using scaled vector transport, we often obtain good convergence properties for the Riemannian conjugate gradient methods.

Scaled vector transport $\mathcal{T}^{(k)}$ satisfies the following inequalities:
\begin{align}~\label{eq:scaled_p1}
\left|\ip{{\grad}f(x_{k+1}),\mathcal{T}^{(k)}_{\alpha_k\eta_k}(\eta_k)}_{x_{k+1}}\right| \leq \left|\ip{{\grad}f(x_{k+1}),\mathcal{T}^{R}_{\alpha_k\eta_k}(\eta_k)}_{x_{k+1}}\right|
\end{align}
and
\begin{align}~\label{eq:scaled_p2}
\norm{\mathcal{T}^{(k)}_{\alpha_k\eta_k}(\eta_k)}_{x_{k+1}} \leq \norm{\eta_{k}}_{x_{k}}.
\end{align}
Now, we would like to verify that inequality \eqref{eq:scaled_p1} holds. From the definition of scaled vector transport \eqref{eq:scaled}, we obtain
\begin{align*}
\left|\ip{{\grad}f(x_{k+1}),\mathcal{T}^{(k)}_{\alpha_k\eta_k}(\eta_k)}_{x_{k+1}}\right|=\left|\ip{{\grad}f(x_{k+1}),s^{(k)}\mathcal{T}^{R}_{\alpha_k\eta_k}(\eta_k)}_{x_{k+1}}\right|,
\end{align*}
where $s^{(k)}$ denotes
\begin{align*}
s^{(k)}:=\min\left\{1,\frac{\|\eta_k\|_{x_k}}{\|\mathcal{T}^R_{\alpha_k\eta_k}(\eta_k)\|_{x_{k+1}}}\right\} \leq 1.
\end{align*}
Therefore, it follows that
\begin{align*}
\left|\ip{{\grad}f(x_{k+1}),\mathcal{T}^{(k)}_{\alpha_k\eta_k}(\eta_k)}_{x_{k+1}}\right|&=s^{(k)}\left|\ip{{\grad}f(x_{k+1}),\mathcal{T}^{R}_{\alpha_k\eta_k}(\eta_k)}_{x_{k+1}}\right| \\
&\leq\left|\ip{{\grad}f(x_{k+1}),\mathcal{T}^{R}_{\alpha_k\eta_k}(\eta_k)}_{x_{k+1}}\right|,
\end{align*}
which leads to \eqref{eq:scaled_p1}. Obviously, \eqref{eq:scaled} implies \eqref{eq:scaled_p2}.

Throughout this paper, we will replace vector transport $\mathcal{T}$ by scaled vector transport $\mathcal{T}^{(k)}$ in \eqref{eq:CGdir}. Therefore, the $(k+1)$-th search direction of the Riemannian conjugate gradient method is determined by
\begin{align}
\label{eq:CGsdir}
\eta_{k+1}=-{\grad}f(x_{k+1})+\beta_{k+1}\mathcal{T}^{(k)}_{\alpha_{k}\eta_{k}}(\eta_{k}).
\end{align}
In \eqref{eq:CGsdir}, $\beta_{k+1}$ is also given by generalizations of the formulas \eqref{eq:EuFR}, \eqref{eq:EuDY}, \eqref{eq:EuPRP}, and \eqref{eq:EuHS}, i.e.,
\begin{align}
\label{eq:FR}
\beta_{k}^\textrm{FR}&=\frac{\|{\grad}f(x_{k})\|^2_{x_{k}}}{\|{\grad}f(x_{k-1})\|^2_{x_{k-1}}}, \\
\label{eq:DY}
\beta_{k}^\textrm{DY}&=\frac{\|{\grad}f(x_{k})\|^2_{x_{k}}}{\ip{{\grad}f(x_{k}),{\cal T}^{(k-1)}_{\alpha_{k-1}\eta_{k-1}}(\eta_{k-1})}_{x_{k}}-\ip{{\grad}f(x_{k-1}),\eta_{k-1}}_{x_{k-1}}}, \\
\label{eq:PRP}
\beta_{k}^\textrm{PRP}&=\frac{\ip{{\grad}f(x_{k}),{\grad}f(x_{k})-{\cal T}^{(k-1)}_{\alpha_{k-1}\eta_{k-1}}({\grad}f(x_{k-1})}_{x_{k}}}{\|{\grad}f(x_{k-1})\|^2_{x_{k-1}}}, \\
\label{eq:HS}
\beta_{k}^\textrm{HS}&=\frac{\ip{{\grad}f(x_{k}),{\grad}f(x_{k})-{\cal T}^{(k-1)}_{\alpha_{k-1}\eta_{k-1}}({\grad}f(x_{k-1})}_{x_{k}}}{\ip{{\grad}f(x_{k}),{\cal T}^{(k-1)}_{\alpha_{k-1}\eta_{k-1}}(\eta_{k-1})}_{x_{k}}-\ip{{\grad}f(x_{k-1}),\eta_{k-1}}_{x_{k-1}}}.
\end{align}
We call these formulas the Fletcher-Reeves, Dai-Yuan, Polak-Ribi\`ere-Polyak and Hestenes-Stiefel formulas, respectively. In the next section, we propose a new choice of $\beta_k$.

In \cite{sato2013}, Sato and Iwai proved that by using the scaled vector transport $\mathcal{T}^{(k)}$ substitute of $\mathcal{T}$ in \eqref{eq:CGdir} and a step size which satisfies the strong Wolfe conditions \eqref{eq:Armijo} and \eqref{eq:sWolfe}, the Fletcher-Reeves type conjugate gradient method defined by \eqref{eq:CGupd}, \eqref{eq:CGsdir}, and \eqref{eq:FR} generates sequences that converge globally. Similarly, in \cite{sato2015}, Sato indicated that if we use scaled vector transport, with a step size satisfying the Wolfe conditions, \eqref{eq:Armijo} and \eqref{eq:Wolfe}, the Dai-Yuan type conjugate gradient method defined by \eqref{eq:CGupd}, \eqref{eq:CGsdir}, and \eqref{eq:DY} generates globally convergent sequences.

\section{Riemannian Hybrid Conjugate Gradient Method and Its Global Convergence Analysis}~\label{sec:main}
\subsection{Proposed hybrid Riemannian conjugate gradient method}
This section describes the Riemannian conjugate gradient descent method using a hybrid $\beta_k$, which exploits the idea described in \cite{dai2001}.

Let $r_k$ be the size of $\beta_k$ with respect to $\beta^\textrm{DY}_k$ defined by \eqref{eq:DY}, namely,
\begin{align}
\label{eq:ratio}
r_k:=\frac{\beta_k}{\beta^\textrm{DY}_k}.
\end{align}
We will prove that, for the method defined by \eqref{eq:CGupd} and \eqref{eq:CGsdir}, the search direction $\eta_k$ is a descent direction at every iteration and the method converges globally if the step size $\alpha_k>0$ satisfies the strong Wolfe conditions \eqref{eq:Armijo} and \eqref{eq:sWolfe}, and the scalar $\beta_k$ is such that
\begin{align}
\label{eq:convcond}
-\sigma \leq r_k \leq 1,
\end{align}
where $\sigma:=(1-c_2)/(1+c_2)>0$ and $c_2$ denotes the constant in the second of the strong Wolfe conditions \eqref{eq:sWolfe}. Furthermore, since the following two choices of $\beta_k$:
\begin{align}
\label{eq:hybrid1}
\beta_k=\max\{0,\min\{\beta^\textrm{DY}_k,\beta^\textrm{HS}_k\}\}
\end{align}
and
\begin{align}
\label{eq:hybrid2}
\beta_k=\max\{-\sigma\beta^\textrm{DY},\min\{\beta^\textrm{DY}_k,\beta^\textrm{HS}_k\}\}
\end{align}
satisfy the condition \eqref{eq:convcond}, we can use either of these hybrid formulas $\beta_k$ defined by \eqref{eq:hybrid1} and \eqref{eq:hybrid2} as the scalar in \eqref{eq:CGsdir}. The above two choices of $\beta_k$ are examples of the hybrid methods in Euclidean space \cite{dai2001}. This implies that our hybrid method is a generalization of the method in \cite{dai2001}. The parameter \eqref{eq:hybrid1} is used in the numerical experiments of \cite[(43)]{hawe2013analysis} and \cite[Table 1]{selvan2012descent}. The hybrid methods using \eqref{eq:hybrid1} and \eqref{eq:hybrid2} combine the good global convergence properties of the Dai-Yuan method \eqref{eq:DY} with the efficient numerical performance of the Hestenes-Stiefel method \eqref{eq:HS}. Now, we note that, in Euclidean space, the hybrid methods using \eqref{eq:hybrid1} and \eqref{eq:hybrid2} converge globally under the Wolfe conditions \eqref{eq:EuWolfe1} and \eqref{eq:EuWolfe2}, whereas, on a Riemannian manifold, the hybrid methods need the strong Wolfe conditions \eqref{eq:Armijo} and \eqref{eq:sWolfe} to converge globally. In Section \ref{sec:exp}, we provide a numerical evaluation showing that the Riemannian conjugate gradient methods with the hybrid $\beta_k$ defined by \eqref{eq:hybrid1} and \eqref{eq:hybrid2} perform better than the Polak-Ribi\`ere-Polyak method.

\subsection{Global convergence analysis}
Zoutendijk's theorem is described on Riemannian manifolds as follows:
\begin{theorem}
[Zoutendijk \cite{ring2012}]
\label{thm:Zoutendjik}
Let $(M,g)$ be a Riemannian manifold and R be a retraction on M.
Let $f:M\rightarrow\real$ be a smooth, bounded below function with the following property: there exists $L>0$ such that
	\begin{align*}
	|D(f \circ R_x)(t\eta)[\eta]-D(f \circ R_x)(0_x)[\eta]| \leq Lt \, \left(\eta \in T_xM,\, \norm{\eta}_x=1,\, x \in M, t \geq 0\right).
	\end{align*}
Suppose that in the line search optimization algorithm \eqref{eq:CGupd}, each step size $\alpha_k>0$ satisfies the strong Wolfe conditions \eqref{eq:Armijo} and \eqref{eq:sWolfe}, and each search direction $\eta_k$ is a descent direction. Then the following series converges:
	\begin{align}
	\label{eq:Zoutendjik}
	\sum_{k=0}^\infty\frac{\ip{{\grad}f(x_k),\eta_k}_{x_k}^2}{\|\eta_k\|^2_{x_k}}<\infty .
	\end{align}
\end{theorem}

The proof of this theorem is along the lines of Zoutendijk's theorem in Euclidean space (see \cite[Theorem 3.3]{ring2012}). Next, we will prove the main convergence theorem.

\begin{theorem}~\label{thm:main}
Let $f:M\rightarrow\real$ be a function satisfying the assumptions of Zoutendijk's theorem. If each $\alpha_k>0$ satisfies the strong Wolfe conditions \eqref{eq:Armijo} and \eqref{eq:sWolfe}, and if $\beta_k$ is such that\footnote{The formulas defined by \eqref{eq:hybrid1} and \eqref{eq:hybrid2} satisfy $-\sigma \leq r_k \leq 1$.} $-\sigma\leq r_k \leq 1$, then any sequence $\{x_k\}$ generated by the Riemannian conjugate gradient method defined by \eqref{eq:CGupd} and \eqref{eq:CGsdir} satisfies
	\begin{align}
	\label{eq:main}
	\liminf_{k \to \infty}\norm{{\grad}f(x_k)}_{x_k}=0.
	\end{align}
\end{theorem}

Let us start with a brief outline of the proof strategy of Theorem \ref{thm:main}, with an emphasis on the main difficulty that has to be overcome in order to generalize the proof in \cite[Theorem 2.3]{dai2001} to manifolds.
The flow of our proof is the same as in \cite{dai2001}.
First, we show that the search direction in each iteration of the hybrid methods is the descent direction.
Therefore, the assumption, "each search direction $\eta_k$ is a descent direction", of Zoutendijk's theorem is satisfied.
Then, assuming that equation \eqref{eq:main} does not hold, the proof is completed by deriving a contradiction.

In general Riemannian manifolds, the inner product of tangent vectors at different points cannot be defined, so the inner product is taken using scaled vector transport. However, the use of scaled vector transport causes a problem that does not occur in Euclidean space. Specifically, the absolute value is required for the inequality in \eqref{eq:abs} when generalizing to the Riemannian manifold.

\begin{proof}[Theorem \ref{thm:main}]
If ${\grad}f(x_{k_0})=0$ for some $k_0$, then \eqref{eq:main} follows.
Thus, it is sufficient to prove \eqref{eq:main} only when ${\grad}f(x_k) \neq 0$ for all $k \geq 0$.

First, we prove that each search direction $\eta_k$ is a descent direction by induction. For $\eta_0=-{\grad}f(x_0)$, it is obvious that $\eta_0$ is a descent direction.

Assume that $\eta_{k-1}$ is a descent direction. Then, we find that
\begin{align}~\label{eq:CGtemp5}
&\ip{{\grad}f(x_k),\eta_k}_{x_k} \nonumber \\
&=\ip{{\grad}f(x_k),-{\grad}f(x_k)+\beta_k{\cal T}^{(k-1)}_{\alpha_{k-1}\eta_{k-1}(\eta_{k-1})}}_{x_k} \nonumber \\
&=-\|{\grad}f(x_k)\|^2_{x_k}+r_k\frac{\|{\grad}f(x_{k})\|^2_{x_{k}}\ip{{\grad}f(x_k),{\cal T}_{\alpha_{k-1}\eta_{k-1}}^{(k-1)}(\eta_{k-1})}_{x_k}}{\ip{{\grad}f(x_{k}),{\cal T}^{(k-1)}_{\alpha_{k-1}\eta_{k-1}}(\eta_{k-1})}_{x_{k}}-\ip{{\grad}f(x_{k-1}),\eta_{k-1}}_{x_{k-1}}} \nonumber \\
\begin{split}
&=\frac{\ip{{\grad}f(x_{k-1}),\eta_{k-1}}_{x_{k-1}}}{\ip{{\grad}f(x_{k}),{\cal T}^{(k-1)}_{\alpha_{k-1}\eta_{k-1}}(\eta_{k-1})}_{x_{k}}-\ip{{\grad}f(x_{k-1}),\eta_{k-1}}_{x_{k-1}}}\|{\grad}f(x_k)\|^2_{x_k} \\
& \qquad +\frac{(r_k-1)\ip{{\grad}f(x_k),{\cal T}^{(k-1)}_{\alpha_{k-1}\eta_{k-1}}(\eta_{k-1})}_{x_{k}}}{\ip{{\grad}f(x_{k}),{\cal T}^{(k-1)}_{\alpha_{k-1}\eta_{k-1}}(\eta_{k-1})}_{x_{k}}-\ip{{\grad}f(x_{k-1}),\eta_{k-1}}_{x_{k-1}}}\|{\grad}f(x_k)\|^2_{x_k}
\end{split}
\end{align}
where the first equation comes from \eqref{eq:CGsdir} and the second equation comes from $\beta_k = r_k\beta_k^\textrm{DY}$ and \eqref{eq:DY}. Accordingly, \eqref{eq:DY} ensures that
\begin{align*}
&\ip{{\rm grad}f(x_k),\eta_k}_{x_k} \\
&=\left\{\ip{{\rm grad}f(x_{k-1}),\eta_{k-1}}_{x_{k-1}}+(r_k-1)\ip{{\rm grad}f(x_{k}),{\cal T}^{(k-1)}_{\alpha_{k-1}\eta_{k-1}}(\eta_{k-1})}_{x_k}\right\}\beta_k^\textrm{DY},
\end{align*}
which, together with \eqref{eq:ratio}, implies that
\begin{align*}
\beta_k&=r_k\beta_k^{\rm DY} \\
&=\frac{r_k\ip{{\grad}f(x_k),\eta_k}_{x_k}}{\ip{{\grad}f(x_{k-1}),\eta_{k-1}}_{x_{k-1}}+(r_k-1)\ip{{\grad}f(x_k),{\cal T}^{(k-1)}_{\alpha_{k-1}\eta_{k-1}}(\eta_{k-1})}_{x_{k}}}.
\end{align*}
Let $l_k$ and $\xi_k$ be
\begin{align}
\label{eq:l}
l_{k}:=&\frac{\ip{{\grad}f(x_k),{\cal T}^{(k-1)}_{\alpha_{k-1}\eta_{k-1}}(\eta_{k-1})}_{x_{k}}}{\ip{{\grad}f(x_{k-1}),\eta_{k-1}}_{x_{k-1}}}, \\
\label{eq:xi}
\xi_k:=&\frac{r_k}{1+(r_k-1)l_{k}}.
\end{align}
Using \eqref{eq:l} and \eqref{eq:xi}, we obtain
\begin{align}~\label{eq:CGtemp4}
\beta_k&=r_k\beta_k^{\rm DY} \nonumber \\
&=\xi_k\frac{\ip{{\grad}f(x_k),\eta_k}_{x_k}}{\ip{{\grad}f(x_{k-1}),\eta_{k-1}}_{x_{k-1}}}.
\end{align}
Furthermore, let $\zeta_k$ be
\begin{align}
\zeta_k:=\frac{1+(r_k-1)l_{k}}{l_k-1}.
\end{align}
Then, \eqref{eq:CGtemp5} guarantees that
\begin{align}\label{eq:CGtemp1}
\ip{{\grad}f(x_k),\eta_k}_{x_k} =\zeta_k\norm{{\grad}f(x_k)}_{x_k}^2.
\end{align}
On the other hand, since $\alpha_k$ satisfies the strong Wolfe conditions, \eqref{eq:sWolfe} implies that
\begin{align*}
\left|\ip{{\grad}f(x_k),\textrm{D}R_{x_{k-1}}(\alpha_{k-1}\eta_{k-1})[\eta_{k-1}]}_{x_k}\right| \leq c_2\left|\ip{{\grad}f(x_{k-1}),\eta_{k-1}}_{x_{k-1}}\right|,
\end{align*}
which, together with \eqref{eq:difft}, \eqref{eq:scaled_p1} and \eqref{eq:l} implies that
\begin{align}~\label{eq:abs}
\begin{split}
|l_{k}|&=\frac{\left|\ip{{\grad}f(x_k),{\cal T}^{(k-1)}_{\alpha_{k-1}\eta_{k-1}}(\eta_{k-1})}_{x_{k}}\right|}{\left|\ip{{\grad}f(x_{k-1}),\eta_{k-1}}_{x_{k-1}}\right|} \\
&\leq\frac{\left|\ip{{\grad}f(x_k),{\cal T}^{R}_{\alpha_{k-1}\eta_{k-1}}(\eta_{k-1})}_{x_{k}}\right|}{\left|\ip{{\grad}f(x_{k-1}),\eta_{k-1}}_{x_{k-1}}\right|} \\
&=\frac{\left|\ip{{\grad}f(x_k),\textrm{D}R_{x_{k-1}}(\alpha_{k-1}\eta_{k-1})[\eta_{k-1}]}_{x_{k}}\right|}{\left|\ip{{\grad}f(x_{k-1}),\eta_{k-1}}_{x_{k-1}}\right|} \leq c_2.
\end{split}
\end{align}
This means $|l_k| \leq c_2 < 1$, which implies $l_k - 1< 0$. Similar to equation (2.18) in \cite{dai2001}, we obtain $1 + (r_k-1)l_k > 0$. Hence,
\begin{align*}
\zeta_k=\frac{1+(r_k-1)l_k}{l_k-1}<0,
\end{align*}
which, together with \eqref{eq:CGtemp1}, implies that $\eta_k$ is a descent direction. Thus, induction shows that each $\eta_k$ is a descent direction.

Finally, we prove \eqref{eq:main} by contradiction. Assume that
\begin{align*}
\liminf_{k \to \infty}\norm{{\grad}f(x_k)}_{x_k}>0.
\end{align*}
Then, noting $\norm{{\grad}f(x_k)}_{x_k} \neq 0$ for all $k$, there exists $\gamma > 0$ such that
\begin{align*}
\|{\grad}f(x_k)\|_{x_k} \geq \gamma >0.
\end{align*}
for all $k$. Since \eqref{eq:CGsdir} means that
\begin{align*}
\eta_k + {\grad}f(x_k) = \beta_k{\cal T}_{\alpha_{k-1}\eta_{k-1}}^{(k-1)}(\eta_{k-1}),
\end{align*}
taking the norms of the above equation and its square, it follows that
\begin{align*}
\|\eta_k\|_{x_k}^2=\beta_k^2\|{\cal T}_{\alpha_{k-1}\eta_{k-1}}^{(k-1)}(\eta_{k-1})\|^2_{x_k}-2\ip{{\grad}f(x_k),\eta_k}_{x_k}-\|{\grad}f(x_k)\|_{x_k}^2.
\end{align*}
Similar to equation (2.21) in \cite{dai2001}, by dividing both sides of the above equation by $\ip{{\grad}f(x_k),\eta_k}_{x_k}^2 \neq 0$, \eqref{eq:CGtemp4} and \eqref{eq:CGtemp1} give,
\begin{align}~\label{eq:CGtemp3}
\begin{split}
\frac{\|\eta_k\|_{x_k}^2}{\ip{{\grad}f(x_k),\eta_k}_{x_k}^2}&=\xi_k^2\frac{\|{\cal T}_{\alpha_{k-1}\eta_{k-1}}^{(k-1)}(\eta_{k-1})\|^2_{x_k}}{\ip{{\grad}f(x_{k-1}),\eta_{k-1}}_{x_{k-1}}^2} \\
& \qquad +\frac{1}{\|{\grad}f(x_k)\|_{x_k}^2}\left\{1-\left(1+\frac{1}{\zeta_k}\right)^2\right\} \textcolor{red}{.}
\end{split}
\end{align}
Similar to equation (2.24) in \cite{dai2001}, we obtain
\begin{align*}
|1 + (r_k - 1)l_k| \geq |r_k|,
\end{align*}
which, together with \eqref{eq:xi}, implies
\begin{align*}
|\xi_k| \leq 1.
\end{align*}
From the above inequality with \eqref{eq:CGtemp3} and \eqref{eq:scaled_p2}, we obtain
\begin{align*}
\frac{\|\eta_k\|_{x_k}^2}{\ip{{\grad}f(x_k),\eta_k}_{x_k}^2} &\leq \frac{\|{\cal T}_{\alpha_{k-1}\eta_{k-1}}^{(k-1)}(\eta_{k-1})\|^2_{x_k}}{\ip{{\grad}f(x_{k-1}),\eta_{k-1}}_{x_{k-1}}^2}+\frac{1}{\|{\grad}f(x_k)\|_{x_k}^2} \\
& \leq \frac{\norm{\eta_{k-1}}_{x_{k-1}}^2}{\ip{{\grad}f(x_{k-1}),\eta_{k-1}}_{x_{k-1}}^2}+\frac{1}{\|{\grad}f(x_k)\|_{x_k}^2}.
\end{align*}
Using the above inequality recursively and noting the hypothesis, $\|{\grad}f(x_k)\|_{x_k} \geq \gamma >0$, and $\|\eta_0\|_{x_0}^2 = \|{\grad}f(x_0)\|_{x_0}^2$, it follows that
\begin{align*}
\frac{\|\eta_k\|_{x_k}^2}{\ip{{\grad}f(x_k),\eta_k}_{x_k}^2} \leq \sum_{i=0}^k\frac{1}{\|{\grad}f(x_i)\|_{x_i}^2} \leq \sum_{i=0}^k\frac{1}{\gamma^2}=\frac{k+1}{\gamma^2}.
\end{align*}
This means
\begin{align*}
\frac{\ip{{\grad}f(x_k),\eta_k}_{x_k}^2}{\|\eta_k\|^2_{x_k}} \geq \frac{\gamma^2}{k+1},
\end{align*}
which indicates
\begin{align*}
\sum_{k=0}^\infty\frac{\ip{{\grad}f(x_k),\eta_k}_{x_k}^2}{\|\eta_k\|_{x_k}^2} \geq \sum_{k=0}^\infty\frac{\gamma^2}{k+1}=\infty .
\end{align*}
This contradicts \eqref{eq:Zoutendjik} in Zoutendijk's theorem and completes the proof.
\qed
\end{proof}

\section{Numerical Experiments}~\label{sec:exp}
This section compares the performances of the existing Riemannian conjugate gradient methods with those of the proposed methods. We solved 7 types of Riemann optimization problems (Problem \ref{pbl:rayleigh}--\ref{pbl:completion}) on several manifolds and objective functions. We solved these problems 10 times with each algorithm, that is, 70 times in total. Then, we calculated a performance profile \cite{dolan2002benchmarking} for each algorithm to show the advantages of our algorithms. Our experiments used the source code of \texttt{pymanopt} (\url{https://github.com/pymanopt}, see \cite{townsend2016pymanopt}). In particular, the Riemannian conjugate gradient method was implemented in \texttt{pymanopt}, so we changed only the parameter $\beta_k$ for the experiments.

\subsection{The Rayleigh-quotient minimization problem on the unit sphere}
Problem \ref{pbl:rayleigh} is the Rayleigh-quotient minimization problem on the unit sphere (see \cite[Chapter 4.6]{absil2008}). The optimal solutions of Problem \ref{pbl:rayleigh} are the unit eigenvectors of $A$ associated with the smallest eigenvalue (see \cite[Chapter 2]{absil2008}).
\begin{problem}~\label{pbl:rayleigh}
For $A \in \mathcal{S}^n_{++}$,
\begin{align*}
\textrm{minimize}& \quad f(x)=x^\top Ax, \\
\textrm{subject to}& \quad x \in \mathbb{S}^{n-1} := \{ x \in \real^{n} : \norm{x} = 1\},
\end{align*}
where $\mathcal{S}^n_{++}$ denotes the set of all symmetric positive-definite matrices.
\end{problem}
In the experiments, we set $n=100$ and generated a matrix $A \in \mathcal{S}^n_{++}$ with randomly chosen elements by using \texttt{sklearn.datasets.make\_spd\_matrix}.

\subsection{Computation of Stability Number}
For an undirected graph $G$, a stable set in $G$ is a set of vertices, which are mutually nonadjacent. We define $S(G)$ as the size of a maximum stable set in $G$. In \cite{motzkin1965maxima}, Motzkin and Straus showed that the computation of the stability number of graphs problem is equivalent to Problem \ref{pbl:stability}. Specifically, the value of the objective function in the global optimal solution of Problem \ref{pbl:stability} is equal to $S(G)^{-1}$. In addition, Yuan, Gu, Lai, and Wen \cite[Section 5.3]{yuan2019global} considered the problem as a Riemannian optimization problem.
\begin{problem}~\label{pbl:stability}
Let $G = (V, E)$ be an undirected graph.
\begin{align*}
\textrm{minimize}& \quad f(x) = \sum_{i=1}^nx_i^4 + 2\sum_{(i,j) \in E}x_i^2x_j^2, \\
\textrm{subject to}& \quad x \in \mathbb{S}^{n-1} := \{ x \in \real^{n} : \norm{x} = 1\},
\end{align*}
where $n=|V|$ and $\norm{\cdot}$ denotes the Euclidean norm.
\end{problem}
In the experiments, we set $n=20$ and generated a graph $G=(V,E)$ randomly by using \texttt{networkx.fast\_gnp\_random\_graph}. Here, we set the probability for edge creationto  $1/4$.

\subsection{The brockett-cost-function minimization problem on a Stiefel manifold}
Problem \ref{pbl:brockett} is the Brockett-cost-function minimization problem on a Stiefel manifold (see \cite[Chapter 4.8]{absil2008}).
\begin{problem}~\label{pbl:brockett}
For $A \in \mathcal{S}^n_{++}$ and $N = \mathrm{diag}(\mu_0, \cdots ,\mu_p)$ $(0 \leq \mu_0\leq \cdots \leq\mu_p)$,
\begin{align*}
\textrm{minimize}& \quad f(X)=\mathrm{tr}(X^\top AXN)\\
\textrm{subject to}& \quad X \in \mathrm{St}(p,n) := \{X \in \real^{n \times p} : X^\top X=I_p\}.
\end{align*}
\end{problem}
In the experiments, we set $p=5$, $n=20$ and $N := \mathrm{diag}(1, \cdots ,p)$ and generated a matrix $A \in \mathcal{S}^n_{++}$ with randomly chosen elements by using \texttt{sklearn.datasets.make\_spd\_matrix}.

\subsection{The closest unit norm column approximation problem}
Problem \ref{pbl:unit-approx} is the closest unit norm column approximation problem, whose implementation is given in \texttt{pymanopt/examples/closest\_unit\_norm\_column\_approximation.py}.
\begin{problem}~\label{pbl:unit-approx}
For $A \in \real^{m \times n}$,
\begin{align*}
\textrm{minimize}& \quad f(X)=\norm{X - A}^2_F \\
\textrm{subject to}& \quad X \in \mathcal{OB}(m,n):=\{X \in \real^{m \times n} : \mathrm{ddiag}(X^\top X)=I_m\},
\end{align*}
where $\norm{\cdot}_F$ denotes the Frobenius norm and $\mathrm{ddiag}(X)$ denotes a diagonal matrix whose diagonal elements are those of $X$.
\end{problem}
In the experiments, we set $m=10$ and $n=1000$ and generated a matrix $A \in \real^{m \times n}$ with randomly chosen elements by using \texttt{numpy.random.randn}.

\subsection{Off-diagonal cost function minimization}
In \cite[Section 3]{absil2006joint}, Absil and Gallivan introduced a cost function on oblique manifolds, which is an off-diagonal cost function written as
\begin{align*}
f(X) := \sum_{i=1}^N\norm{X^\top C_iX - \mathrm{ddiag}(X^\top C_iX)}^2_F,
\end{align*}
where $C_i$ $(i=1,2, \cdots ,N)$ are symmetric matrices. Problem \ref{pbl:off-diag} is one of minimizing the off-diagonal cost function on an oblique manifold.
\begin{problem}~\label{pbl:off-diag}
For $C_i \in \mathcal{S}^{n}$ $(i=1, \cdots ,N)$,
\begin{align*}
\textrm{minimize}& \quad f(X)=\sum_{i=1}^N\norm{X^\top C_iX - \mathrm{ddiag}(X^\top C_iX)}^2_F \\
\textrm{subject to}& \quad X \in \mathcal{OB}(n,p):=\{X \in \real^{n \times p} : \mathrm{ddiag}(X^\mathrm{T}X)=I_p\},
\end{align*}
where $\mathcal{S}^n$ denotes the set of all symmetric matrices.
\end{problem}
In the experiments, we set $N=5$, $n=10$ and $p=5$ and generated 5 matrices $B_i \in \real^{n \times n}$ $(i=1,2, \cdots ,5)$ with randomly chosen elements by using \texttt{numpy.random.randn}. We set symmetric matrices $C_i \in \mathcal{S}^n$ as $C_i := (B_i + B_i^\top) / 2$ $(i=1,2, \cdots ,5)$.

\subsection{The low-rank matrix approximation problem}
Problem \ref{pbl:low-approx} is the low-rank matrix approximation problem whose implementation is given in \texttt{pymanopt/examples/low\_rank\_matrix\_approximation.py}.
\begin{problem}~\label{pbl:low-approx}
For $A \in \real^{m \times n}$,
\begin{align*}
\textrm{minimize}& \quad f(X)=\norm{X - A}^2_F, \\
\textrm{subject to}& \quad X \in M_k := \{ X \in \real^{m \times n} : \mathrm{rank}(X) = k\}.
\end{align*}
\end{problem}
In the experiments, we set $m=100$, $n=80$ and $k=4$ and generated a matrix $A \in \real^{m \times n}$ with randomly chosen elements by using \texttt{numpy.random.randn}.

\subsection{The robust matrix completion problem}
Problem \ref{pbl:completion} is the robust matrix completion problem, discussed by Vandereycken \cite[Section 1.1 (1.5)]{vandereycken2013low}.
\begin{problem}~\label{pbl:completion}
For $A \in \real^{m \times n}$, and a subset $\Omega$ of the complete set of entries $\{1, \cdots ,m\}\times\{1, \cdots ,n\}$,
\begin{align*}
\textrm{minimize}& \quad f(X)=\norm{P_\Omega(X - A)}^2_F, \\
\textrm{subject to}& \quad X \in M_k := \{ X \in \real^{m \times n} : \mathrm{rank}(X) = k\},
\end{align*}
where
\begin{align*}
P_\Omega : \real^{m \times n} \rightarrow \real^{m \times n}, X_{ij} \mapsto
\begin{cases}
X_{ij} & (i,j) \in \Omega \\
0 & (i,j) \not\in \Omega
\end{cases}.
\end{align*}
\end{problem}
In the experiments, we set $m=10$, and $m=8$ and $k=4$, and $\Omega$ contained each pair $(i,j) \in \{1, \cdots ,m\}\times\{1, \cdots ,n\}$ with probability $1/2$. Moreover, we used a matrix $A \in \real^{m \times n}$ that was generated with randomly chosen elements by using \texttt{numpy.random.randn}.

We used line search algorithms for the strong Wolfe conditions \eqref{eq:Armijo} and \eqref{eq:sWolfe} with $c_1=0.0001$ and $c_2=0.9$. We determined that a sequence had converged to an optimal solution if the stopping condition,
\begin{align*}
\norm{{\grad}f(x_k)}_{x_k}<10^{-6}
\end{align*}
was satisfied.

The experiments used a MacBook Air (2017) with a 1.8 GHz Intel Core i5, 8 GB 1600 MHz DDR3 memory, and macOS Mojave version 10.14.5 operating system. The algorithms were written in Python 3.7.6 with the NumPy 1.17.3 package and the Matplotlib 3.1.1 package. We modified the strong Wolfe line search provided as scipy.optimize.line\_search in the SciPy package, to compute the step size in \eqref{eq:CGupd}.

For comparison, we chose two Riemannian conjugate gradient methods, i.e., the Dai-Yuan method \eqref{eq:DY} and the Polak-Ribi\`ere-Polyak method \eqref{eq:PRP}. Below, we call the hybrid methods using \eqref{eq:hybrid1} and \eqref{eq:hybrid2}, Hybrid1 and Hybrid2, respectively.

Table \ref{tb:summary1} and \ref{tb:summary2} summarize the results such as the average and median values of the above 70 experiments. In particular, Table \ref{tb:summary1} shows summary statistics for the number of iterations and Table \ref{tb:summary2} shows those for the elapsed time. From Table \ref{tb:summary1} and \ref{tb:summary2}, we can see that the hybrid methods converge to optimal solutions in fewer iterations and in less time than the DY and PRP methods.
\begin{table}[htbp]
\centering
\caption{Summary statistics on the iteration of 70 experiments of the Riemannian optimization problems. \label{tb:summary1}}
\begin{tabular}{c||c|c|c|c}
\hline\noalign{\smallskip}
 & DY & PRP & Hybrid1 & Hybrid2 \\
\noalign{\smallskip}\hline\noalign{\smallskip}
mean & 1438.7 & 399.9 & 212.2 & 235.0 \\
std & 1765.9 & 513.8 & 217.0 & 212.0 \\
min & 46 & 21 & 20 & 20 \\
median & 570.5 & 161 & 129 & 135 \\
max & 7061 & 2037 & 952 & 803 \\
\noalign{\smallskip}\hline
\end{tabular}
\end{table}
\begin{table}[htbp]
\centering
\caption{Summary statistics on the elapsed time of 70 experiments of the Riemannian optimization problems. \label{tb:summary2}}
\begin{tabular}{c||c|c|c|c}
\hline\noalign{\smallskip}
 & DY & PRP & Hybrid1 & Hybrid2 \\
\noalign{\smallskip}\hline\noalign{\smallskip}
mean & 18.70 & 4.39 & 2.56 & 2.91 \\
std & 26.22 & 4.53 & 2.13 & 2.46 \\
min & 0.43 & 0.14 & 0.15 & 0.14 \\
median & 10.34 & 2.74 & 2.31 & 2.44 \\
max & 147.46 & 20.95 & 10.74 & 11.42 \\
\noalign{\smallskip}\hline
\end{tabular}
\end{table}

Then, we calculate the performance profiles \cite{dolan2002benchmarking}. The performance profile $P_s:\real \rightarrow [0,1]$ is defined as follows: let $\mathcal{P}$ and $\mathcal{S}$ be the set of problems and solvers, respectively. For each $p \in \mathcal{P}$ and $s \in \mathcal{S}$, we define $t := \textrm{(computing time required to solve problem }p\textrm{ by solver }s)$. We define the performance ratio $r_{p,s}$ as
\begin{align*}
r_{p,s}:= \dfrac{t_{p,s}}{\min_{s^\prime \in \mathcal{S}}t_{p,s^\prime}}.
\end{align*}
Next, we define the performance profile, for all $\tau \in \real$, as
\begin{align*}
P_s(\tau) := \frac{\mathrm{size}\{p \in \mathcal{P} : r_{p,s} \leq \tau \}}{\mathrm{size}\mathcal{P}},
\end{align*}
where $\mathrm{size}A$ denotes the number of elements of a set $A$.

\begin{figure}[htbp]
 \centering
 \includegraphics[scale=0.35]{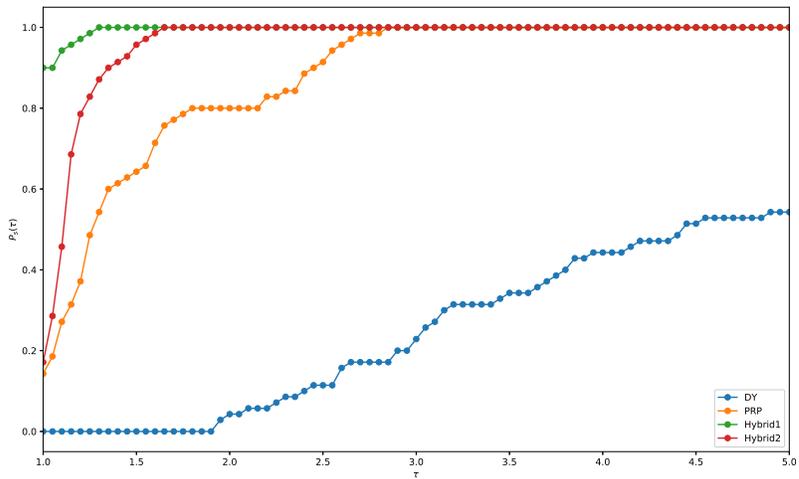}
\caption{Performance profile versus number of iterations. \label{fig1}}
\end{figure}

\begin{figure}[htbp]
 \centering
 \includegraphics[scale=0.35]{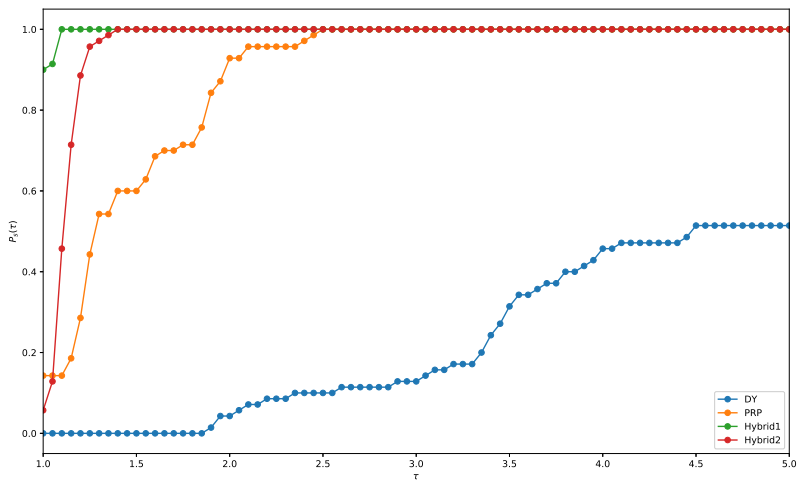}
\caption{Performance profile versus elapsed time. \label{fig2}}
\end{figure}

Figure \ref{fig1} plots the performance profile of each algorithm versus the number of iterations. It shows that the hybrid methods have much higher performance than the DY method. Moreover, the hybrid methods outperform the PRP method. Also, it can be seen that Hybrid1 is superior to Hybrid2. Figure \ref{fig2} plots the performance profiles of each algorithm versus the elapsed time. We can see that the hybrid methods are superior to both DY and PRP. In particular, they perform much better than the DY method. In addition, Hybrid1 is again superior to Hybrid2.

\section{Conclusion and Future Work}~\label{sec:conclude}
This paper presented hybrid Riemannian conjugate gradient methods and showed their global convergence properties.
It compared them numerically with the existing Riemannian conjugate gradient methods on several Riemannian optimization problems.
The results of the numerical experiments demonstrated the efficiency of the hybrid methods.

Various hybrid conjugate methods have been proposed for Euclidean space, such as
\begin{align*}
\beta_k=\max\{ 0, \min\{\beta_k^\textrm{PRP}, \beta_k^\textrm{FR}\} \}.
\end{align*}
The hybrid conjugate methods in Euclidean space are summarized in \cite{hager2006}. We will present more hybrid methods and convergence analyses in a future paper.

\section{Acknowledgment}
We are sincerely grateful to the editor and the anonymous reviewer for helping us improve the original manuscript. This work was supported by JSPS KAKENHI Grant Number JP18K11184.

\bibliographystyle{abbrv}
\bibliography{biblib}

\begin{thebibliography}{10}

\bibitem{absil2006joint}
P.-A. Absil and K.~A. Gallivan.
\newblock Joint diagonalization on the oblique manifold for independent
  component analysis.
\newblock In {\em 2006 IEEE International Conference on Acoustics Speech and
  Signal Processing Proceedings}, volume~5, pages V--V, 2006.

\bibitem{absil2008}
P.-A. Absil, R.~Mahony, and R.~Sepulchre.
\newblock {\em Optimization {A}lgorithms on {M}atrix {M}anifolds}.
\newblock Princeton University Press, 2008.

\bibitem{baali1985}
M.~Al-Baali.
\newblock Descent property and global convergence of the {F}letcher-{R}eeves
  method with inexact line search.
\newblock {\em IMA Journal of Numerical Analysis}, 5(1):121--124, 1985.

\bibitem{dai1999}
Y.-H. Dai and Y.~Yuan.
\newblock A nonlinear conjugate gradient method with a strong global
  convergence property.
\newblock {\em SIAM Journal on optimization}, 10(1):177--182, 1999.

\bibitem{dai2001}
Y.-H. Dai and Y.~Yuan.
\newblock An efficient hybrid conjugate gradient method for unconstrained
  optimization.
\newblock {\em Annals of Operations Research}, 103(1-4):33--47, 2001.

\bibitem{dolan2002benchmarking}
E.~D. Dolan and J.~J. Mor{\'e}.
\newblock Benchmarking optimization software with performance profiles.
\newblock {\em Mathematical programming}, 91(2):201--213, 2002.

\bibitem{fletcher1964}
R.~Fletcher and C.~M. Reeves.
\newblock Function minimization by conjugate gradients.
\newblock {\em The computer journal}, 7(2):149--154, 1964.

\bibitem{hager2006}
W.~W. Hager and H.~Zhang.
\newblock A survey of nonlinear conjugate gradient methods.
\newblock {\em Pacific journal of Optimization}, 2(1):35--58, 2006.

\bibitem{hawe2013analysis}
S.~Hawe, M.~Kleinsteuber, and K.~Diepold.
\newblock Analysis operator learning and its application to image
  reconstruction.
\newblock {\em IEEE Transactions on Image Processing}, 22(6):2138--2150, 2013.

\bibitem{hestenes1952}
M.~R. Hestenes and E.~Stiefel.
\newblock {\em Methods of conjugate gradients for solving linear systems}.
\newblock NBS Washington, DC, 1952.

\bibitem{hu1991}
Y.~Hu and C.~Storey.
\newblock Global convergence result for conjugate gradient methods.
\newblock {\em Journal of Optimization Theory and Applications},
  71(2):399--405, 1991.

\bibitem{motzkin1965maxima}
T.~S. Motzkin and E.~G. Straus.
\newblock Maxima for graphs and a new proof of a theorem of turán.
\newblock {\em Canadian Journal of Mathematics}, 17:533–540, 1965.

\bibitem{polak1969}
E.~Polak and G.~Ribi{\`e}re.
\newblock Note sur la convergence de m{\'e}thodes de directions conjugu{\'e}es.
\newblock {\em ESAIM: Mathematical Modelling and Numerical
  Analysis-Mod{\'e}lisation Math{\'e}matique et Analyse Num{\'e}rique},
  3(R1):35--43, 1969.

\bibitem{ring2012}
W.~Ring and B.~Wirth.
\newblock Optimization methods on {R}iemannian manifolds and their application
  to shape space.
\newblock {\em SIAM Journal on Optimization}, 22(2):596--627, 2012.

\bibitem{sato2015}
H.~Sato.
\newblock A {D}ai-{Y}uan-type {R}iemannian conjugate gradient method with the
  weak {W}olfe conditions.
\newblock {\em Computational Optimization and Applications}, 64(1):101--118,
  2016.

\bibitem{sato2013}
H.~Sato and T.~Iwai.
\newblock A new, globally convergent {R}iemannian conjugate gradient method.
\newblock {\em Optimization}, 64(4):1011--1031, 2015.

\bibitem{selvan2012descent}
S.~E. Selvan, U.~Amato, K.~A. Gallivan, C.~Qi, M.~F. Carfora, M.~Larobina, and
  B.~Alfano.
\newblock Descent algorithms on oblique manifold for source-adaptive ica
  contrast.
\newblock {\em IEEE Transactions on Neural Networks and Learning Systems},
  23(12):1930--1947, 2012.

\bibitem{smith1994}
S.~T. Smith.
\newblock Optimization techniques on {R}iemannian manifolds.
\newblock {\em Fields Institute Communications}, 3(3):113--135, 1994.

\bibitem{touati1990}
D.~Touati-Ahmed and C.~Storey.
\newblock Efficient hybrid conjugate gradient techniques.
\newblock {\em Journal of Optimization Theory and Applications},
  64(2):379--397, 1990.

\bibitem{townsend2016pymanopt}
J.~Townsend, N.~Koep, and S.~Weichwald.
\newblock Pymanopt: A python toolbox for optimization on manifolds using
  automatic differentiation.
\newblock {\em The Journal of Machine Learning Research}, 17(1):4755--4759,
  2016.

\bibitem{vandereycken2013low}
B.~Vandereycken.
\newblock Low-rank matrix completion by {R}iemannian optimization.
\newblock {\em SIAM Journal on Optimization}, 23(2):1214--1236, 2013.

\bibitem{wolfe1969}
P.~Wolfe.
\newblock Convergence conditions for ascent methods.
\newblock {\em SIAM Review}, 11(2):226--235, 1969.

\bibitem{wolfe1971}
P.~Wolfe.
\newblock Convergence conditions for ascent methods. ii: Some corrections.
\newblock {\em SIAM Review}, 13(2):185--188, 1971.

\bibitem{yuan2019global}
H.~Yuan, X.~Gu, R.~Lai, and Z.~Wen.
\newblock Global optimization with orthogonality constraints via stochastic
  diffusion on manifold.
\newblock {\em Journal of Scientific Computing}, 80(2):1139--1170, 2019.

\end{thebibliography}

\end{document}